\documentclass[12pt]{article}
\usepackage{amsmath,amsfonts,calrsfs,amssymb,color}

\newtheorem{Theorem}{Theorem}[section]
\newtheorem{Definition}[Theorem]{Definition}

\newtheorem{Lemma}[Theorem]{Lemma}

\makeatletter
\newif\ifmsbmloaded@

\def\loadmsbm{\msbmloaded@true
 \font\tenmsb=msbm10 scaled 1\@ptsize00
 \font\sevenmsb=msbm7 scaled 1\@ptsize00
 \font\fivemsb=msbm5 scaled 1\@ptsize00
 \alloc@8\fam\chardef\sixt@@n\msbfam
 \textfont\msbfam=\tenmsb
 \scriptfont\msbfam=\sevenmsb
 \scriptscriptfont\msbfam=\fivemsb
 }
\def\nonmatherr@#1{\errmessage%
{LateX error: \string#1\space allowed only in math mode}}
\def\Bbb{\relax\ifmmode\expandafter\Bbb@\else
 \expandafter\nonmatherr@\expandafter\Bbb\fi}
\def\Bbb@#1{{\Bbb@@{#1}}}
\def\Bbb@@#1{\fam\msbfam\relax#1}

\@addtoreset{equation}{section}

\makeatother
\loadmsbm

\def\R{\mathbb R}
\def\N{\mathbb N}

\def\E{\mathbb E}
\def\P{\mathbb P}

\def\ds{\displaystyle}

\title{Existence and uniqueness of nonnegative solutions to the stochastic porous media equation}
   \author{Viorel Barbu \thanks{Supported by the CEEX Project 05 of
Romanian Minister of Research.},\\
University Al. I. Cuza \\and \\
Institute of Mathematics ``Octav Mayer'', Iasi, Romania
,\\
Giuseppe Da Prato \thanks{Supported by the research program ``Equazioni di
Kolmogorov'' from the Italian
``Ministero della Ricerca Scientifica e Tecnologica''},\\
 Scuola Normale Superiore
di Pisa, Italy\\
 and \\
Michael R\"ockner \thanks{Supported  by the
DFG-Research Group 399, the SFB-701, the
BIBOS-Research Center, the INTAS project 99-559,
the RFBR project 04--01--00748,
the Russian--Japanese Grant 05-01-02941-JF,
the DFG Grant 436 RUS 113/343/0(R).} \\
Faculty of Mathematics, University of Bielefeld, Germany\\
and\\
Department of Mathematics and Statistics, Purdue University,\\  U. S. A.}

\begin{document}
 \maketitle

{\bf Abstract}. One proves that the stochastic porous media equation in $3$-D has a unique nonnegative solution for nonnegative
initial data in $H^{-1}(\mathcal O)$ if the nonlinearity is monotone and has polynomial growth.



\section{Introduction} 
Let $\mathcal O$ be an open bounded domain of $\R^n$ with smooth boundary  $\partial\mathcal O.$
We consider the  linear operator $\Delta$ in $L^2(\mathcal O)$ defined on $H^2(\mathcal O)\cap H_0^1(\mathcal O).$
It is well known that $-\Delta$ is self-adjoint positive and anti-compact. So, there exists
 a complete orthonormal system $\{e_k\}$ in $L^2(\mathcal O)$ of eigenfunctions  of $-\Delta$. 
In fact we have $e_k\in \cap_{p\ge 1}
L^p(\mathcal O)$ for all $k\in \N$. We denote by
$\{\lambda_k\}$ the corresponding sequence of eigenvalues,
$$
\Delta e_k=-\lambda_k e_k,\quad k\in \N.
$$
 We shall consider a cylindrical Wiener process in $L^2(\mathcal O)$  of the following form
$$
W(t)=\sum_{k=1}^\infty \beta_k(t)e_k,\quad t\ge 0,
$$
where    $\{\beta_k\}$ is
a sequence of mutually independent   standard Brownian
motions on a filtered probability space
$(\Omega,\mathcal F,\{\mathcal F_t\}_{t\ge 0},\P)$.  
To be  more specific, we shall assume that $1\le n\le 3$.

In this work we consider the stochastic partial differential equation,
\begin{equation}
\label{e1.1}
\left\{\begin{array}{l}
dX(t)-\Delta\beta(X(t))dt=XdW(t),\quad t\ge 0,\\
\\
\beta(X(t))=0,\quad\mbox{\rm on}\; \partial \mathcal O,\quad t\ge 0,\\
\\
X(0,x)=x.
\end{array}\right.
\end{equation}
 Here $\beta$ is a continuous, differentiable,  monotonically increasing function on $\R$ which satisfies the
following conditions,
\begin{equation}
\label{e1.2}
\left\{\begin{array}{l}
|\beta'(r)|\le \alpha_1|r|^{m-1}+\alpha_2,\quad \forall\;r\in \R,\\
\\
\ds j(r)\colon=\int_0^r\beta(s)ds\ge\alpha_3|r|^{m+1}+\alpha_4r^2,\quad \forall\;r\in \R,
\end{array}\right.
\end{equation}
where $\alpha_i>0,\;i=1,2,3,4$ and $1\le m$.
We note that since $\beta$ is increasing, the mean value theorem implies that
\begin{equation}
\label{e1.3}
r\beta(r)\ge j(r),\quad r\ge 0.
\end{equation} 

Equation \eqref{e1.1} with additive noise was recently studied in \cite{BBDR},\cite{BDR},\cite{DRRW}, 
\cite{DR1},\cite{DR2}, see also \cite{BD}. In particular, in
\cite{DRRW} was given an
existence result under similar conditions on $\beta$. 
Here we consider a multiplicative noise (of a special form, but it would be possible to consider a
more general noise $f(X)dW(t)$ with $f(0)=0$), which is needed in
order to ensure   positivity of  solutions.

As was shown in \cite{PR} existence and uniqueness of solutions follow
by the general results in \cite{PR} (see also \cite{RRW} for generalizations). In this paper we present an alternative proof,
based on the  Yosida approximation of
$-\Delta\beta$,  and prove the positivity of solutions for nonnegative initial data $x$.

As in deterministic case the Sobolev space $H^{-1}(\mathcal O)$ is   natural  for studying equation
\eqref{e1.1}. Equation \eqref{e1.1} can be written in the abstract form  
\begin{equation}
\label{e1.4}
\left\{\begin{array}{l}
dX(t)+AX(t)=\sigma(X(t))dW(t),\quad t\ge 0,\\
\\
X(0)=x,
\end{array}\right.
\end{equation}
where the operator $A\colon D(A)\subset H^{-1}(\mathcal O)\to H^{-1}(\mathcal O)$ is defined by
\begin{equation}
\label{e1.5}
\left\{\begin{array}{l}
Ax=-\Delta\beta(x),\quad x\in D(A),\\
\\
D(A)=\{x\in H^{-1}(\mathcal O)\cap L^1(\mathcal O):\; \beta(x)\in  H^{1}_0(\mathcal O)   \},
\end{array}\right.
\end{equation}
and where
\begin{equation}
\label{e1.6}
\sigma(X)dW(t)=\sum_{k=1}^\infty\mu_k Xe_kd\beta_k(t),\quad X\in H^{-1}(\mathcal O).
\end{equation}
To give a rigorous sense to this noise term we first note that since $n\le 3$, by Sobolev embedding it follows that 
\begin{equation}
\label{e1.7}
\sup_{k\in \N}\frac1{\lambda_k}\;|e_k|_\infty<\infty.
\end{equation}
Furthermore, troughout this paper we shall assume that
\begin{equation}
\label{e1.8}
\sum_{k=1}^\infty\mu^2_k \lambda^2_k=\colon C<\infty.
\end{equation}
  \eqref{e1.8}  implies for some constant $c_1>0$
\begin{equation}
\label{e1.9}
 \sum_{k=1}^\infty\mu^2_k|xe_k|^2_{-1}\le  c_1\sum_{k=1}^\infty\mu^2_k\;\lambda^2_k\;|x|^2_{-1}\le c_1C|x|^2_{-1},\quad \forall\; x\in
H^{-1}(\mathcal O),
\end{equation}
because $|xe_k|^2_{-1}\le c_1\lambda^2_k|x|^2_{-1}$ by an elementary calculation, since $n\le 3$ and due to
\eqref{e1.7}.

Defining
\begin{equation}
\label{e1.10}
 \sigma(x)h:=\sum_{k=1}^\infty\mu_k(h,e_k) xe_k,\quad x\in H^{-1}(\mathcal O),\;h\in L^2(\mathcal O),
\end{equation}
we obtain by \eqref{e1.9} that $ \sigma(x)\in L_2(L^2(\mathcal O),H^{-1}(\mathcal O))$. Considering $(\beta_k)_{k\in \N}$
as a cylindrical Wiener process on $L^2(\mathcal O)$, it follows that \eqref{e1.6} is well defined. Note that since
$\sigma$ is linear we have that $x\to \sigma(x)$ is Lipschitz from $H^{-1}(\mathcal O)$ to $L_2(L^2(\mathcal
O),H^{-1}(\mathcal O))$ (in particular \cite{KR}, \cite{PR}, \cite{RRW} really apply).
\bigskip

The plan of the paper  is the following:   main results are stated in \S 2  and  proofs are given in \S 3.

The following notations will be used throughout in the following.

\begin{enumerate}
\item[(i)] $H^1_0(\mathcal O),H^2(\mathcal O)$ are standard Sobolev spaces on $\mathcal O$ endowed with their usual norms
denoted by $|\cdot|_{H^1_0(\mathcal O)}$ and $|\cdot|_{H^2(\mathcal O)}$ respectively.

\item[(ii)] $H$ is the space $H^{-1}(\mathcal O)$ (the dual of $H^1_0(\mathcal O)$) endowed with the norm
$$
|x|_H=|x|_{-1}=|-\Delta^{-1}x|_{H^1_0(\mathcal O)}.
$$
(Here $(-\Delta)^{-1}x=y$ is the solution to Dirichlet problem $-\Delta y=x$ in $\mathcal O$, $y\in H^1_0(\mathcal O)$).
The scalar product in $H$ is
$$
\langle x,z  \rangle_{-1}=\int_\mathcal O(-\Delta)^{-1}\,x \,z \,d\xi,\quad \forall\; x,z\in H^1_0(\mathcal O).
$$

\item[(iii)] The scalar product and the norm in $L^2(\mathcal O)$ will be denoted by $(\cdot,\cdot)$
and $|\cdot|_2$, respectively and the norm in $L^p(\mathcal O),\;1\le p\le \infty$  by
$|\cdot|_p$.

\item[(iv)] For two Hilbert spaces $H_1$, $H_2$ the space of Hilbert-Schmidt operators from $H_1$ to $H_2$ is denoted by  $L_2(H_1,H_2)$.
\end{enumerate}

\section{The main result}
To begin with let us define the solution concept we shall work with.
Formally, a solution to \eqref{e1.1} (equivalently \eqref{e1.4}) might be an $H$--valued 
continuous adapted process
such that $X,AX\in C_W([0,T];L^2(\Omega;H))$ and
\begin{equation}
\label{e2.1}
X(t)=x-\int_0^tAX(s)ds+\int_0^t \sigma(X(s))dW(s),\quad t\in [0,T].
\end{equation}
(By $C_W([0,T];L^2(\Omega;H))$ we mean the  Banach space of all the processes $X$ in $(\Omega,\mathcal F,\P)$
with values in $H$ which are adapted and mean square continuous, endowed with the norm
$$
\|X\|_{C_W([0,T];L^2(\Omega;H))}^2:=\sup_{t\in[0,T]}\E|X(t)|^2_H.
$$
Spaces $L^p_W([0,T];L^2(\Omega;H))$,
$p\in [1,\infty]$, are defined similarly.)

However,  such a concept of solution might fail to exist for equation \eqref{e1.1} and so we shall confine to a
weaker  one inspired by \cite{DRRW} and \cite{KR}.
\begin{Definition}
\label{d2.1}
An $H$-valued continuous $\mathcal F_t$-adapted process $X$ is called a solution to \eqref{e1.1} on $[0,T]$ if $X\in
L^{m+1} (\Omega\times(0,T)\times \mathcal O )$ and
\begin{equation}
\label{e2.2}
\begin{array}{lll}
(X(t),e_j)&=&\ds(x,e_j)+\int_0^t\int_\mathcal O\beta(X(s))\Delta e_j d\xi ds\\
\\
&&\ds +\sum_{k=1}^\infty\mu_k\int_0^t(X(s)e_k,e_j)d\beta_k(s),\quad \forall\;j\in \N,\;t\in [0,T].
\end{array}
\end{equation}
\end{Definition}
Taking into account that $-\Delta e_j =\lambda_j e_j $ in $\mathcal O$ we may equivalently write \eqref{e2.2} as
follows
$$
\begin{array}{lll}
\langle X(t),e_j\rangle_{-1}&=&\ds\langle x,e_j\rangle_{-1}- \int_0^t\int_\mathcal O\beta(X(s))e_j d\xi ds\\
\\
&&\ds +\sum_{k=1}^\infty\mu_k\int_0^t\langle X(s)e_k,e_j\rangle_{-1}d\beta_k(s),\quad \forall\;j\in \N,
\end{array}
$$
i.e.
$$
d\langle X(t),e_j\rangle_{-1}+ (\beta(X(t)),e_j)dt
=\sum_{k=1}^\infty\mu_k\langle X(s)e_k,e_j\rangle_{-1}d\beta_k(s).
$$
Recalling \eqref{e1.6} we see that
$$
\sum_{k=1}^\infty\mu_k(X(t)e_k,e_j)d\beta_k(t)=(\sigma(X(t))dW(t),e_j),\quad j\in \N.
$$
We   also note that since by assumption \eqref{e1.2}, $\beta(X)\in L^{\frac{m+1}m}((0,T)\times \Omega\times
\mathcal O)$, the integral arising in the right hand side of \eqref{e2.2} makes sense because $e_j\in 
C^\infty(\overline{\mathcal O})$ for all $j\in \N$. 
Of course, one might derive a vector valued version of Definition  \ref{d2.1} as in \cite{DRRW}.
Now we are ready to formulate the main results.
\begin{Theorem}
\label{t2.2}
Assume that \eqref{e1.2} and \eqref{e1.8} hold. Then for each $x\in H^{-1}(\mathcal O)$ there is a unique solution $X$
to 
\eqref{e1.1}.
Moreover, if $x\in L^p(\mathcal O)$ is nonnegative a.e. on $\mathcal O$ where $p\geq \max\{m+1,4\}$ is a natural
number then $X\in
L^\infty_W(0,T;L^p(\Omega;L^p(\mathcal O)))$ and 
$X\ge 0$ a.e. on
$(0,\infty)\times \mathcal O$, $\P$-a.s. If $x\in H^{-1}(\mathcal O)$ is such that $x\ge 0$, i.e. $x$ is a positive
measure,  then $X(t)\ge 0$ for all $t\ge 0$, $\P$-a.s.
\end{Theorem} 
The positivity of the solution $X$ to \eqref{e1.1} will be proven below by choosing an appropriate Lyapunov function.

\section{Proof of Theorem \ref{t2.2}}

We mention that in our estimates in the sequel constants may change from line to line though we do not express this in our notation.

We recall that the operator $A$, defined by \eqref{e1.5}, is maximal monotone in $H$ (see e.g. \cite{B}). Then we consider the Yosida approximation
$$
A_\varepsilon(x)=\frac1\varepsilon\;(x-J_ \varepsilon (x))=A(1+\varepsilon A)^{-1}(x),\quad \varepsilon>0,\;x\in H,
$$
where $J_ \varepsilon (x)=(1+\varepsilon A)^{-1}(x).$
The operator $A_\varepsilon$ is monotone and Lipschitzian on $H$. Then, by \eqref{e1.9}
it follows by standard existence theory for stochastic equations in the Hilbert spaces (see e.g. \cite{DPZ1}) that
 the approximating equation
\begin{equation}
\label{e3.1}
\left\{\begin{array}{l}
dX_\varepsilon(t)+A_\varepsilon X_\varepsilon(t)dt=\sigma(X_\varepsilon(t))dW(t),\quad t\ge 0,\\
\\
X_\varepsilon(0)=x,
\end{array}\right.
\end{equation}
has a unique solution $X_\varepsilon\in C_W([0,T];L^2(\Omega;H))$ such that $X_\varepsilon\in C([0,T];H), \P$-a.s. with 
$A_\varepsilon X_\varepsilon\in C_W([0,T];L^2(\Omega;H))$.\bigskip

  By It\^o's formula we have
\begin{equation}
\label{e3.2}
\begin{array}{l}
\ds \frac12\;d|X_\varepsilon(t)|^2_{-1}+\langle A_\varepsilon X_\varepsilon(t),
X_\varepsilon(t)  \rangle_{-1}dt
\\
\\
\ds=\langle \sigma (X_\varepsilon(t))dW(t),
X_\varepsilon(t)  \rangle_{-1}+\frac12\;\sum_{k=1}^\infty\mu_k^2|X_\varepsilon(t) e_k|^2_{-1}dt.
\end{array}
\end{equation}
This yields (see    \eqref{e1.9})
$$
\begin{array}{l}
\ds\frac12\;\E|X_\varepsilon(t)|^2_{-1}+\E\int_0^t\langle A_\varepsilon X_\varepsilon(s),
X_\varepsilon(s)  \rangle_{-1}ds\\
\\
\ds \le\frac12\;|x|^2_{-1}+C\E\int_0^t|X_\varepsilon(s)|^2_{-1}ds
\end{array}
$$
and therefore
\begin{equation}
\label{e3.3}
\frac12\;\E|X_\varepsilon(t)|^2_{-1}+\E\int_0^t\langle A_\varepsilon X_\varepsilon(s),
X_\varepsilon(s)  \rangle_{-1}ds\le C|x|^2_{-1},\quad\forall\;\varepsilon>0.
\end{equation}
We set $Y_\varepsilon(t)=J_\varepsilon(X_\varepsilon(t))$ (see \eqref{e3.1}).  Then
\begin{equation}
\label{e3.4}
\begin{array}{l}
\ds\frac12\;\E|X_\varepsilon(t)|^2_{-1}+\E\int_0^t\int_\mathcal Oj  (Y_\varepsilon(s))ds d\xi\\
\\
\ds+\frac1\varepsilon\;\E\int_0^t|X_\varepsilon(s)-Y_\varepsilon(s)|^2_{-1}ds\le
C|x|^2_{-1},\quad\forall\;\varepsilon>0.
\end{array}
\end{equation}
(Here we have used the equality
$$
\langle A_\varepsilon x,x  \rangle_{-1}=\langle A  J_\varepsilon x,J_\varepsilon x  \rangle_{-1}+
\frac1\varepsilon\;|x-J_\varepsilon (x)|^2_{-1},
$$
and \eqref{e1.3}.)
 
Now we fix $X\in C_W([0,T];L^2(\Omega,H))$ and we consider the equation
\begin{equation}
\label{e3.5}
\left\{\begin{array}{l}
d\widetilde{X}_\varepsilon(t)+A_\varepsilon \widetilde{X}_\varepsilon(t)dt=\sigma(X(t))dW(t),\quad t\ge 0,\\
\\
\widetilde{X}_\varepsilon(0)=x.
\end{array}\right.
\end{equation}
Equivalently,
\begin{equation}
\label{e3.6}
\left\{\begin{array}{l}
d\widetilde{X}_\varepsilon(t)-\Delta\beta(\widetilde{Y}_\varepsilon(t))dt=\sigma(X(t))dW(t),\quad t\ge 0,\\
\\
\widetilde{X}_\varepsilon(0)=x,
\end{array}\right.
\end{equation}
where
$$
\widetilde{Y}_\varepsilon=(1+\varepsilon A)^{-1}\widetilde{X}_\varepsilon.
$$
 On the other hand, for equation  \eqref{e3.5} we have the same estimates as for  \eqref{e3.1}. In fact by It\^o's formula
 we get (see  \eqref{e3.4})
\begin{equation}
\label{e3.7}
\begin{array}{lll}
\ds\E|\widetilde{X}_\varepsilon(t)|^2_{-1} +\E\int_0^t\int_\mathcal Oj  (\widetilde{Y}_\varepsilon(s))ds d\xi &+&\ds
\frac1\varepsilon\;\E\int_0^t|\widetilde{X}_\varepsilon(s)-\widetilde{Y}_\varepsilon(s)|^2_{-1}ds
\\
&\le&\ds C|x|^2_{-1}
+C\E\int_0^t|X(s)|^2_{-1},
\end{array}
\end{equation}
(where we have used \eqref{e1.9} to estimate $\E\int_0^t\|\sigma(X (s))\|_{L_2(L^2(\mathcal O);H)}^2ds$).
By virtue of assumption \eqref{e1.2} this  implies that
$$
\E\int_0^T\int_\mathcal O 
 |\beta(\widetilde{Y}_\varepsilon(s))|^{\frac{m+1}m}ds d\xi\le C(|x|^2_{-1}+1),\quad \varepsilon>0,
$$
(because $|\beta(r)|\le \tilde\alpha_1|r|^m+\tilde\alpha_2,\;\tilde\alpha_1\ge 0$),
and so along a subsequence, we have
\begin{equation}
\label{e3.8}
\beta(\widetilde{Y}_\varepsilon)\to \eta\quad\mbox{\rm weakly in}\;L^{\frac{m+1}{m}}((0,T)\times \Omega\times
\mathcal O).
\end{equation}
On the other hand, we have by \eqref{e3.6} that for   $t\in [0,T]$
$$
\langle \widetilde{X}_\varepsilon(t),e\rangle_{-1}+\int_0^t\int_\mathcal O\beta(\widetilde{Y}_\varepsilon(s)) e d\xi ds
=\langle x,e\rangle _{-1}+\int_0^t\langle \sigma(X(s))dW(s),e\rangle _{-1}ds, 
$$
for all $e\in L^{m+1}(\mathcal O)$.  We note that  by \eqref{e3.7} there exists $X^*\in L^2_W([0,T];L^2(\Omega;H))$
such that
\begin{equation}
\label{e3.10}
\widetilde{X}_\varepsilon\to X^*\quad\mbox{\rm weakly in }\;L^2_W([0,T];L^2(\Omega;H))
\end{equation}
and by \eqref{e3.7} and \eqref{e1.2} we obtain that also
\begin{equation}
\label{e3.10bis}
\widetilde{Y}_\varepsilon\to X^*\quad\mbox{\rm weakly in }\;L^2_W([0,T];L^2(\Omega;H))\cap
 L^{m+1}(\Omega\times(0,T)\times\mathcal
O).
\end{equation}
Hence along a subsequence $\varepsilon \to 0$
$$
\E\langle \widetilde{X}_\varepsilon(t),e\rangle _{-1}\to \E\langle X^*(t),e\rangle _{-1}\quad \mbox{\rm weakly
in}\;L^2(0,T).
$$
Then letting $\varepsilon$ tend to $0$ we get for a.e.  $t\in [0,T]$
\begin{equation}
\label{e3.9}
\langle X^*(t),e\rangle _{-1} =\langle
x,e\rangle_{-1}-\int_0^t\int_\mathcal O\eta(s) ed\xi ds+\int_0^t\langle \sigma(X(s))dW(s),e\rangle _{-1}ds.
\end{equation}

Taking into account \eqref{e3.10}-\eqref{e3.10bis}, to conclude the proof of existence it suffices to show that
\begin{equation}
\label{e3.11}
\eta(t,\xi,\omega)=\beta(X^*(t,\xi,\omega))\quad\mbox{\rm a.e.}\;(\omega,t,\xi)\in
\Omega\times(0,T)\times\mathcal O.
 \end{equation}
Indeed, in such a case we may take in \eqref{e3.9} $e=\Delta e_j$ for $j\in \N$.

To this end we consider the operator
$$
F\colon L^{m}(\Omega\times(0,T)\times 
\mathcal O)\to L^{\frac{m}{m+1}}(\Omega\times(0,T)\times 
\mathcal O)=(L^{m}(\Omega\times(0,T)\times 
\mathcal O))',
$$
defined by
$$
(Fx)(t,\xi,\omega)=\beta(x(t,\xi,\omega))\quad\mbox{\rm a.e.}\;(\omega,t,\xi)\in \Omega\times(0,T)\times\mathcal O.
$$
This operator is maximal monotone and more precisely, it is the subgradient of the convex function
$\Phi:L^{m+1}(\Omega\times(0,T)\times  \mathcal O)\to \R$ defined as,
$$
\Phi(x)=\E\int_0^T\int_\mathcal Oj(x(t,\xi))dt d\xi.
$$
For each $Z\in L^{m+1}(\Omega\times(0,T)\times
\mathcal O)$ we have
$$
\Phi(\widetilde{Y}_\varepsilon)-\Phi(Z)\le \E\int_0^T\int_\mathcal O
\beta(\widetilde{Y}_\varepsilon(t,\xi))(\widetilde{Y}_\varepsilon(t,\xi)-
Z(t,\xi)) dt d\xi
$$
Letting $\varepsilon$ tend to $0$ we have by \eqref{e3.8}, \eqref{e3.10}, \eqref{e3.10bis} and by the weak lower semicontinuity of
$\Phi$
$$
\Phi(X^*)-\Phi(Z)\le \liminf_{\varepsilon\to 0}\E\int_0^T\int_\mathcal O
\beta(\widetilde{Y}_\varepsilon(t,\xi))\widetilde{Y}_\varepsilon(t,\xi) dt d\xi
-\E\int_0^T\int_\mathcal O\eta Z dt d\xi.
$$
To prove \eqref{e3.11} by the uniqueness of the subgradient it suffices to show that
\begin{equation}
\label{e3.12}
\liminf_{\varepsilon\to 0}\E\int_0^T\int_\mathcal O
\beta(\widetilde{Y}_\varepsilon(t,\xi))\widetilde{Y}_\varepsilon(t,\xi) dt d\xi\le 
\E\int_0^T\int_\mathcal O\eta X^* dt d\xi.
\end{equation}
To this end we come back to equation \eqref{e3.6} and note that by It\^o's  formula we have
$$
\begin{array}{l}
\ds\frac12\;\E|\widetilde{X}_\varepsilon(t)|^2_{-1}+\E\int_0^t\int_\mathcal O\beta(\widetilde{Y}_\varepsilon(s))
\widetilde{X}_\varepsilon(s)dsd\xi
\\
\\
\ds=\frac12\;|x|^2_{-1}+\frac12\;\sum_{k=1}^\infty\E\int_0^t\mu_k^2|X(s)e_k|^2_{-1}ds.
\end{array}
$$
Equivalently,
\begin{equation}
\label{e3.12'}
\begin{array}{l}
\ds\frac12\;\E|\widetilde{X}_\varepsilon(t)|^2_{-1}+\E\int_0^t\int_\mathcal O\beta(\widetilde{Y}_\varepsilon(s))
\widetilde{Y}_\varepsilon(s)dsd\xi\\
\\
\ds+
\E\int_0^t\int_\mathcal O\beta(\widetilde{Y}_\varepsilon(s))
(\widetilde{X}_\varepsilon(s)-\widetilde{Y}_\varepsilon(s))dsd\xi
\\
\\
\ds=\frac12\;|x|^2_{-1}+\frac12\;\sum_{k=1}^\infty\E\int_0^t\mu_k^2|X(s)e_k|^2_{-1}ds.
\end{array}
\end{equation}
By \eqref{e3.10}-\eqref{e3.10bis} we have
$$
\int_\mathcal O\beta(\widetilde{Y}_\varepsilon(s))
(\widetilde{X}_\varepsilon(s)-\widetilde{Y}_\varepsilon(s))d\xi=
\langle  A_\varepsilon \widetilde{X}_\varepsilon(s),\widetilde{X}_\varepsilon
(s)-J_\varepsilon(\widetilde{X}_\varepsilon(s))
 \rangle_{-1}
=\varepsilon|A_\varepsilon \widetilde{X}_\varepsilon(s)|^2_{-1}.
$$
Fix $\varphi\in L^\infty(0,T)$, $\varphi\ge 0$. Then $\varphi X^*\in L^2_W(0,T;L^2(\Omega;H))$. Thus by \eqref{e3.10}-\eqref{e3.10bis}
$$
\begin{array}{l}
\ds\E\int_0^T\varphi(t)|X^*(t)|^2_{-1}\;dt=\lim_{\varepsilon\to 0}\E\int_0^T\langle X^*(t),X_\varepsilon(t)  \rangle_{-1}
\varphi(t)dt
\\
\\
\ds \le\left( \E\int_0^T\varphi(t)|X^*(t)|^2_{-1}\;dt\right)^{1/2}
\liminf_{\varepsilon\to 0}\left( \E\int_0^T\varphi(t)|X_\varepsilon(t)|^2_{-1}\;dt\right)^{1/2}.
\end{array}
$$
Hence simplifying we obtain
$$
\E\int_0^T\varphi(t)|X^*(t)|^2_{-1}\;dt\le \liminf_{\varepsilon\to 0} \E\int_0^T\varphi(t)|X_\varepsilon(t)|^2_{-1}\;dt.
$$ 
Hence \eqref{e3.12'}, Fatou's Lemma (see also \eqref{e1.3}) and the arbitrariness of $\varphi$ implies that for a.e. $t\in [0,T]$
we obtain that
\begin{equation}
\label{e3.13}
\begin{array}{l}
\ds\liminf_{\varepsilon\to 0}\E\int_0^t\int_\mathcal O
\beta(\widetilde{Y}_\varepsilon(s)\widetilde{Y}_\varepsilon(s) ds d\xi+\frac12\;
\E|X^*(t)|^2_{-1}\\
\\
\ds\le 
\frac12\;|x|^2_{-1}+\frac12\;\sum_{k=1}^\infty\E\int_0^t\mu_k^2|X(s)e_k|^2_{-1}ds.
\end{array}
\end{equation}
On the other hand, by \eqref{e3.9} we see via It\^o's formula (applied to the  right hand side of \eqref{e3.9}, since the left hand side might
not be continuous in $t$)  that for all $j\in \N$ and a.e. $t\in
[0,T]$,
$$
\begin{array}{l}
\ds \frac12\;\E|\langle X^*(t),e_j   \rangle_{-1}|^2+\E\int_0^t\langle\eta_s,e_j   \rangle\langle X^*(s),e_j 
\rangle_{-1} ds
\\
\\
\ds=\frac12\;\langle x,e_j   \rangle^2_{-1}+\frac12\;\E\sum_{k=1}^\infty\mu_k^2\int_0^t\langle X(s)e_k,e_j  \rangle^2
ds
\end{array}
$$
and dividing by $ |e_j|^2_{-1}$ and  summing over $j$ we obtain
\begin{equation}
\label{e3.14}
\begin{array}{l}
\ds\frac12\;\E|X^*(t)|^2_{-1}+\E\int_0^t\int_\mathcal O\eta(s)X^*(s) dsd\xi
\\
\\
\ds=\frac12\;|x|^2_{-1}+\frac12\;\sum_{k=1}^\infty\mu_k^2\E\int_0^t|X(s)e_k|^2_{-1}ds.
\end{array}
\end{equation}
We note that the integral in the left hand side makes sense since by \eqref{e3.4},
$X^*\in L^{m+1}((0,T)\times \Omega\times\mathcal O)$ while $\eta\in L^{\frac{m+1}m}((0,T)\times \Omega\times
\mathcal O)$.

Comparing \eqref{e3.13} and \eqref{e3.14} we infer that
$$
\liminf_{\varepsilon\to 0}\E\int_0^T\int_\mathcal O
\beta(\widetilde{Y}_\varepsilon(t))\widetilde{Y}_\varepsilon(t) dt d\xi\le 
\E\int_0^T\int_\mathcal O\eta(t) X^*(t) dt d\xi,
$$
as claimed.   A 
formal problem arises, however, because $X^*(t)$
 as constructed before might not be $H$-continuous.  However,  arguing as in \cite{KR}, \cite{PR} we may replace it by an
$H$-continuous version defined by
$$
\widetilde{X}^*(t)=x+\int_0^t\Delta\eta(s)ds+ \int_0^t\sigma(X(s))dW(s).
$$
It follows that $X^*=\widetilde{X}^*$ a.e.  and that    $\widetilde{X}^*$
is also  an $\mathcal F_t$-adapted process. Moreover, the It\^o formula  from (\cite[Theorem I-3-2]{KR}) holds.
Hence $\widetilde{X}^*\in C_W([0,T];L^2(\Omega;H))\cap L^{m+1}((0,T)\times \Omega\times \mathcal O) $
is a solution (in the sense of Definition 2.1) to  
\begin{equation}
\label{e3.15}
\left\{\begin{array}{l}
dX^*+AX^*dt=\sigma(X) dW\\
\\
X^*(0)=x.
\end{array}\right.
\end{equation}\bigskip

\bigskip
{\bf Uniqueness}. Let $X^*_1,X^*_2$ be two solutions to equation  \eqref{e1.1} for $X=X_i,\;i=1,2$. We have (see 
\eqref{e2.2})
$$
d\langle X^*_1-X^*_2,e_j\rangle _{-1}+\int_\mathcal O(\beta(X^*_1)-\beta(X^*_2))e_j d\xi dt=
\sum_{k=1}^\infty\mu_k\langle (X_1-X_2)e_k,e_j\rangle _{-1}d\beta_k.
$$
By It\^o's  formula we obtain
$$
\begin{array}{l}
\ds\frac12\;\E|\langle X^*_1(t)-X^*_2(t),e_j\rangle _{-1}|^2\\
\\
\ds+\E\int_0^t( \beta(X^*_1(s))-\beta(X^*_2(s)),e_j)\langle X^*_1(s)-X^*_2(s),e_j\rangle _{-1}ds
\\
\\
\ds=\frac12\;\E\int_0^t\sum_{k=1}^\infty\mu^2_k\langle (X_1(s)-X_2(s))e_k,e_j\rangle ^2_{-1}ds
\end{array}
$$
Dividing by $|e_j|^2_{-1}$ and 
summing over $j$ we see that
$$
\begin{array}{l}
\ds\frac12\;\E|X^*_1(t)-X^*_2(t)|_{-1}^2+\E\int_0^t(\beta(X^*_1)-\beta(X^*_2),X^*_1(s)-X^*_2(s)\;ds
\\
\\
\ds=\frac12\;\E\int_0^t\sum_{j,k=1}^\infty\mu^2_k\langle (X_1(s)-X_2(s))e_k,|e_j|^{-1}_{-1}\,e_j\rangle ^2_{-1}ds.
\end{array}
$$
Hence (see \eqref{e1.9})
\begin{equation}
\label{e3.16}
\E|X^*_1(t)-X^*_2(t)|_{-1}^2\le CE\int_0^t|X_1(s)-X_2(s)|^2_{-1}ds,\quad\forall\;t\in[0,T]
\end{equation}
Now we shall use the latter inequality to prove existence of a unique solution
$$
X\in C_W([0,T];L^2(\Omega;H))\cap L^{m+1}((0,T)\times \Omega\times\mathcal O)
$$
to equation \eqref{e1.1}. Indeed the operator $X\to X^*$ is a contraction on the space 
$C_W([0,T];L^2(\Omega;H))$ if $T$ is sufficiently small and so, we have existence (and uniqueness) for $T>0$
small. By  a  standard unique continuation argument it follows existence and uniqueness on an arbitrary interval
$[0,T]$.

 \bigskip

{\bf Positivity}. We shall assume now that  $x\in L^p(\mathcal O)$, where $p\ge \max\{m+1,4\}$,   and   $x(\xi)\ge 0$ a.e. 
in $\mathcal O$. We shall prove that
\begin{equation}
\label{e3.17}
X\ge 0 \quad\mbox{\rm a.e. in }\;(0,T)\times \mathcal O\times \Omega.
\end{equation}
We shall first  assume in addition  that $\beta$ is strictly monotone, i.e.
\begin{equation}
\label{e3.18}
(\beta(r)-\beta(\bar r))(r-\bar r)\ge \alpha(r-\bar r)^2,\quad\forall\;r,\bar r\in \R,
\end{equation}
where $\alpha>0$. Below we shall use the following lemma.
\begin{Lemma}
\label{l3.1}
Let $y\in D(A)$ and $g:\R\to \R$ Lipschitz and increasing. Then
$$
\langle\nabla\beta(y),\nabla g(y)\rangle_{\R^n} \ge 0,\quad\mbox{\it a.e. on }\;\mathcal O.
$$
\end{Lemma}
{\bf Proof}. First note that by definition of $D(A)$ we have that $y,\beta(y)\in H^1_0(\mathcal O)$.
Using a Dirac sequence we can find mollifiers $g_k\in C^1(\R)$, $g'_k\ge 0,\;k\in \N,$ such that
$$
\nabla g(y)=\lim_{k\to \infty}g'_k(y)\nabla  y\quad\mbox{\rm in}\;L^2(\mathcal O).
$$
So, it suffices to prove that
$$
\langle\nabla\beta(y),  \nabla g(y)\rangle_{\R^n} \ge 0,\quad\mbox{\rm a.e. on}\;\mathcal O.
$$
But
$$
\langle\nabla\beta(y),  \nabla y\rangle_{\R^n}=\langle\nabla\beta(y),  \nabla\beta^{-1}\beta(y)\rangle_{\R^n}.
$$
Since $\beta$ is strictly monotone, $\beta^{-1}$ is Lipschitz, so applying the above mollifier argument with $\beta^{-1}$
replacing $g$, we prove the assertion. $\Box$\bigskip

We shall use the approximating equation \eqref{e3.1} whose solution   $X_\varepsilon$
is weakly convergent to $X$ in $L^2_W(\Omega;L^2(0,T;H))$. Namely, we have
for $Y_\varepsilon(t):=J_\varepsilon(X_\varepsilon(t))$, $t\ge 0$,
\begin{equation}
\label{e3.19}
dX_\varepsilon(t)-\Delta\beta(Y_\varepsilon (t))dt=\sigma(X_\varepsilon(t))dW(t),\quad t\ge 0.
\end{equation}
We note that equation \eqref{e3.1} can be equivalently written as
\begin{equation}
\label{e3.20}
\left\{\begin{array}{l}
\ds dX_\varepsilon(t)+\frac1\varepsilon\; X_\varepsilon(t) dt=\frac1\varepsilon\; J_\varepsilon(X_\varepsilon(t)) dt+
\sigma(X_\varepsilon(t))dW(t),\quad t\ge
0,\\
\\
X_\varepsilon(0)=x,
\end{array}\right.
\end{equation}
Fix $x\in H$ and set
$$
y=J_\varepsilon(x)=(1-\varepsilon\Delta\beta)^{-1}x,
$$
i.e.
\begin{equation}
\label{e3.21}
y-\varepsilon\Delta\beta(y)=x
\end{equation}
Then $y\in D(A)$. Since $\beta$ is strictly monotone, $\beta^{-1}$ is Lipschitz. Therefore, since
$\beta(y)\in H^1_0(\mathcal O)$, also $y\in H^1_0(\mathcal O)$. Now assume 
$x\in L^p(\mathcal O)$. 
By multiplying both sides of \eqref{e3.21} by $\frac{y^{p-1}}{1+\lambda y^{p-2}}$ and
integrating over $\mathcal O$ we get by Lemma \ref{l3.1}
$$
\int_\mathcal O\frac{y^p}{1+\lambda |y|^{p-2}}\;d\xi\le\int_\mathcal O\frac{y^{p-1}x}{1+\lambda |y|^{p-2}}\;d\xi.
$$
Then, letting $\lambda\to 0$ we find the estimate
\begin{equation}
\label{e3.22}
|y|_p^p\le \int_\mathcal O y^{p-1}xd\xi\le |y|^{p-1}_p\;|x|_p.
\end{equation}
Hence
\begin{equation}
\label{e3.23}
|J_\varepsilon(x)|_{p}\le |x|_{p},\quad\forall\;x\in L^p(\mathcal O),
\end{equation}
and therefore,
$$
|A_\varepsilon(x)|_{p}=\frac1\varepsilon\;|x-J_\varepsilon(x)|_p
\le \frac2\varepsilon\;|x|_p,\quad\forall\;x\in L^p(\mathcal O).
$$
\eqref{e3.21} and \eqref{e3.23} imply that $J_\varepsilon$ is continuous from $L^p(\mathcal O)$ into itself.
\begin{Lemma}
\label{l3.2}
For each $x\in L^2(\mathcal O)$ equation   \eqref{e3.20} has a unique solution
 $X_\varepsilon\in C_W([0,T];L^2(\Omega;L^2(\mathcal O))).$
\end{Lemma}
{\bf Proof}. Let us first prove that $J_\varepsilon=(1-\varepsilon \Delta\beta)^{-1}$
is  Lipschitz continuous in $L^2(\mathcal O)$.
Indeed, by the equation
$$
J_\varepsilon(x)-\varepsilon\Delta\beta(J_\varepsilon(x))=x,\quad\mbox{\rm in}\;\mathcal O,
$$
(taking into account that $\beta(J_\varepsilon(x))\in H^1_0(\mathcal O)$) we have for $x,\bar x\in 
L^2(\mathcal O)$
$$
\begin{array}{l}
\ds \int_\mathcal O(J_\varepsilon(x)-J_\varepsilon(\bar x))(\beta(J_\varepsilon(x))-\beta(J_\varepsilon(\bar x)))d\xi
\\
\\
\ds +\varepsilon\int_\mathcal O|\nabla(\beta(J_\varepsilon(x))-\beta(J_\varepsilon(\bar x))|^2d\xi
\le \int_\mathcal O(x-\bar x)(\beta(J_\varepsilon(x))-\beta(J_\varepsilon(\bar x)))d\xi.
\end{array}
$$
This yields, recalling \eqref{e3.18}
$$
\alpha |J_\varepsilon(x)-J_\varepsilon(\bar x)|^2_2+\varepsilon
|\beta(J_\varepsilon(x))-\beta(J_\varepsilon(\bar x))|^2_{H^1_0(\mathcal O)}\le 
|x-\bar x|_2\; |\beta(J_\varepsilon(x))-\beta(J_\varepsilon(\bar x))|_2.
$$
On the other hand, by the Poincar\'e inequality there exists $C>0$ such that
$$
 |\beta(J_\varepsilon(x))-\beta(J_\varepsilon(\bar x))|^2_2\le C|\beta(J_\varepsilon(x))-\beta(J_\varepsilon(\bar x))|^2_{H^1_0(\mathcal O)}.
$$
Therefore
$$
\begin{array}{l}
\ds\alpha |J_\varepsilon(x)-J_\varepsilon(\bar x))|^2_2+\frac\varepsilon2\;
|\beta(J_\varepsilon(x))-\beta(J_\varepsilon(\bar x))|^2_{H^1_0(\mathcal
O)}+\frac\varepsilon{2C}\;|\beta(J_\varepsilon(x))-\beta(J_\varepsilon(\bar x))|^2_2\\
\\
\ds\le \frac{C}{2\varepsilon}\; |x-\bar x|^2_2\;
+\frac\varepsilon{2C}\;|\beta(J_\varepsilon(x))-\beta(J_\varepsilon(\bar x))|^2_2,
\end{array}
$$
and consequently
$$
\alpha |J_\varepsilon(x)-J_\varepsilon(\bar x))|^2_2+\frac\varepsilon2\;
|\beta(J_\varepsilon(x))-\beta(J_\varepsilon(\bar x))|^2_{H^1_0(\mathcal
O)}\le\frac{C}{2\varepsilon}\; |x-\bar x|_2.
$$
So, $J_\varepsilon$
is  Lipschitz continuous in $L^2(\mathcal O)$  as claimed. Consequently
$A_\varepsilon=\frac1\varepsilon\;(1-J_\varepsilon)$ is  Lipschitz continuous in $L^2(\mathcal O)$ as well.
Moreover, since
$$
\begin{array}{l}
\ds\|\sigma(x)\|_{L_2(L^2(\mathcal O),L^2(\mathcal O))}\le \sum_{k=1}^\infty\mu_k^2|xe_k|^2_2 \le 
\sum_{k=1}^\infty\mu_k^2|e_k|^2_{L^\infty(\mathcal O)}\;|x|^2_2\le
C_1\sum_{k=1}^\infty\mu_k^2\lambda_k^2\;|x|^2_2
\end{array}
$$
we infer by standard existence theory for stochastic PDEs that for each $x\in L^2(\mathcal O)$ equation
\eqref{e3.20} has a unique solution in $X_\varepsilon\in C_W([0,T];L^2(\Omega;L^2(\mathcal O)))$ (see e.g.
\cite{DPZ1}). $\Box$ \bigskip

For $R>0$ define
$$
K_R:=\{X\in L^\infty_W(0,T;L^p(\Omega\times \mathcal O)):\;e^{-4\alpha t}\E|X(t)|^p_p\le R^p\quad\mbox{\rm for a.e.}\;t\in [0,T]\}
$$
 \begin{Lemma}
\label{l3.3}
Let $T>0$ and $x\in L^p(\mathcal O)$. Then for the solution $X_\varepsilon$ of \eqref{e3.1} (or equivalently
\eqref{e3.20}) we have $X_\varepsilon\in L^\infty_W(0,T;L^p(\Omega\times\mathcal O))$
and $X_\varepsilon$ is bounded in $L^\infty_W(0,T;L^p(\Omega\times\mathcal O))$
\end{Lemma}
{\bf Proof}. Obviously, $K_R$ is a closed subset of $L^\infty_W(0,T;L^p(\Omega\times\mathcal O)).$
Since by \eqref{e3.20} $X_\varepsilon$ is a fixed point of the map 
$$
X\stackrel{F}{\mapsto} 
 e^{-\frac{t}{\varepsilon}}x+\frac1\varepsilon\;\int_0^te^{-\frac{(t-s)}{\varepsilon}} J_\varepsilon(X
(s))ds+\int_0^te^{-\frac{(t-s)}{\varepsilon}}  \sigma(X (s))dW(s),\quad t\in [0,T],
$$
obtained by iteration in $C_W(0,T;L^2(\Omega\times \mathcal O))$, it suffices to prove that this map
leaves $K_R$ invariant for $R$ large enough. But for $X\in K_R$ we have by \eqref{e3.23} for $t\ge 0$
$$
\begin{array}{l}
\ds\left(e^{-p\alpha t}\E\left|e^{-\frac{t}{\varepsilon}}x+\frac1\varepsilon\;\int_0^te^{-\frac{(t-s)}{\varepsilon}}
J_\varepsilon(X (s))ds  \right|_p^p   \right)^{1/p}
\\
\\
\ds\le e^{-\alpha t}e^{-\frac{t}{\varepsilon}}|x|_p+e^{-\alpha t}\left(\E\left(\int_0^t
\frac1\varepsilon\;e^{-\frac{(t-s)}{\varepsilon}}
|J_\varepsilon(X (s))|_pds \right)^p\right)^{1/p}\\
\\
\ds\le e^{-(\frac1\varepsilon+\alpha) t} |x|_p+e^{-\alpha t}\Bigg(\int_0^t\cdots\int_0^t
\frac1\varepsilon\;e^{-\frac{(t-s_1)}{\varepsilon}}e^{\alpha s_1}\cdots \frac1\varepsilon\;e^{-\frac{(t-s_p)}{\varepsilon}}e^{\alpha s_p}\\
\\
\hspace{30mm}\times e^{-\alpha s_1}(\E|(|X(s_1)|_p^p)^{1/p}\cdots e^{-\alpha s_p}(\E|(|X(s_p)|_p^p)^{1/p}\;ds_1\cdots ds_p\Bigg)^{1/p}\\
\\
\ds \le e^{-(\frac1\varepsilon+\alpha) t} |x|_p+e^{-\alpha t}R\int_0^t\frac1\varepsilon\;e^{-\frac{(t-s)}{\varepsilon}}e^{\alpha s}ds\\
\\
\le 
e^{-(\frac{1}{\varepsilon}+\alpha)t}|x|_p+\frac{R}{1+\alpha\varepsilon}.
\end{array}
$$
Now we set
$$
Y(t)=\int_0^te^{-\frac{(t-s)}{\varepsilon}}X(s))dW(s),\quad t\ge 0.
$$
 Then
$$
\left\{\begin{array}{l}
\ds dY(t)+\frac1\varepsilon\;Y(t)dt=\sigma(X(t))dW(t),\quad t\ge 0,\\
\\
Y(0)=0.
\end{array}\right. 
$$
Let $\lambda>0$. Applying It\^o's formula to the function
$$
\Psi_\lambda(y):=\frac1p\;|(1+\lambda A_0)^{-1}y|^p_p,\quad y\in L^p(\mathcal O),
$$
(see the beginning of the proof of the next lemma for a detailed justification) we obtain via H\"older's inequality that
$$
\begin{array}{l}
\ds \E[\Psi_\lambda(Y(t))]+\frac1\varepsilon\;\E\int_0^t\int_\mathcal O|(1+\lambda A_0)^{-1}Y(s)|^pd\xi\,ds
\\
\\
\ds=\frac{p-1}2\;\sum_{k=1}^\infty\mu_k^2\;\E\int_0^t\int_\mathcal O|(1+\lambda A_0)^{-1}Y(s)|^{p-2}\\
\\
\ds\hspace{10mm} \times |(1+\lambda A_0)^{-1}(X(s)e_k)|^2d\xi\,ds\\
\\
\ds\le C\E\int_0^t|(1+\lambda A_0)^{-1}Y(s)|^2_p\;|X(s)|^2_pds\\
\\
\ds\le \frac1{2\varepsilon}\;\E\int_0^t|(1+\lambda A_0)^{-1}Y(s)|^p_p ds+\frac{9C^2\varepsilon}{8}\;
\E\int_0^t |X(s)|^p_pds\\
\\
\ds\le \frac1{2\varepsilon}\;\E\int_0^t|(1+\lambda A_0)^{-1}Y(s)|^p_p ds+\frac{9C^2\varepsilon(e^{4\alpha t}-1) }{32\alpha}\;R^p.
\end{array}
$$
Then letting $\lambda\to \infty$, we see by Fatou's lemma
that for a.e. $t\in [0,T]$ we have for $C_1$ independent of $\varepsilon$
$$
e^{-4\alpha t}\E|Y(t)|_p^p\le   \frac{C_1\varepsilon}\alpha  R^p ,\quad \forall\;t\in [0,T].
$$
This means that for $\alpha$ large enough and 
  $R>2|x|_p$  the map leaves $K_R$
invariant as claimed. 

\begin{Lemma}
\label{l3.4}
For $x\in L^p(\mathcal O))$ we have
$$
X_\varepsilon\to X\quad\mbox{\it strongly in}\;L^\infty_W(0,T;L^2(\Omega;H)),
$$
$$
X_\varepsilon\to X\quad\mbox{\it weakly in}\;L^\infty_W(0,T;L^p(\Omega;L^p(\mathcal O))),
$$
where $X$ is the solution to \eqref{e1.1}.
\end{Lemma}
{\bf Proof}. By \eqref{e3.4} and Lemma \ref{l3.3} we know that $\{X_\varepsilon\}$ is bounded in
$$
L^2_W(0,T;L^2(\Omega;H))\cap L^\infty_W(0,T;L^p(\Omega;L^p(\mathcal O)))
$$
Subtracting equations \eqref{e1.1} and \eqref{e3.1} we get via It\^o's formula and because $\beta$ is increasing that
$$
\begin{array}{l}
\ds\frac12\;\E|X_\varepsilon(t)-X(t)|^2_{-1}+\E\int_0^t\int_\mathcal O(\beta((1+\varepsilon A)^{-1}X)-\beta(X))
(X_\varepsilon-X)ds d\xi
\\
\\
\ds\le c\E\int_0^t|X_\varepsilon(s)-X(s)|^2_{-1}ds,
\end{array}
$$
and by Gronwall's lemma we obtain
\begin{equation}
\label{e3.24bis}
\E|X_\varepsilon(t)-X(t)|^2_{-1}\le
C\E\int_0^1\int_\mathcal O(\beta((1+\varepsilon A)^{-1}X)-\beta(X))
(X_\varepsilon-X)ds \,d\xi.
\end{equation}
On the other hand, it follows by \eqref{e3.23} that
$$
\int_{\Omega\times [0,T]\times\mathcal O}|(1+\varepsilon A)^{-1}X|^p \P(d\omega)
 \,dt\, d\xi\le
\int_{\Omega\times [0,T]\times\mathcal O}|X|^p \P(d\omega)
 \,dt \,d\xi,
$$
while for $\varepsilon\to 0$
$$
(1+\varepsilon A)^{-1}X\to X\quad\mbox{\rm in }\;L^1(\mathcal O)
$$
for $(\omega,t)\in \Omega\times [0,T]$ (which is a consequence of the fact that the  operator $A$ is $m$-accretive in $L^1(\mathcal O)$,
cfr. \cite{Aubin}).
Hence (at least along a subsequence)
$$
(1+\varepsilon A)^{-1}X\to X\quad\mbox{\rm a.e. on }\;\Omega\times [0,T]\times\mathcal O.
$$
Hence
$$
(1+\varepsilon A)^{-1}X\to X\quad\mbox{\rm weakly in }\;L^p(\Omega\times [0,T]\times\mathcal O)
$$
as $\varepsilon\to 0$ and according to the above inequality this implies that for $\varepsilon\to 0$,
$|(1+\varepsilon A)^{-1}X|_{L^p}\to |X|_{L^p}$. Hence since $L^p(\Omega\times [0,T]\times\mathcal O)$ is uniformly convex,
$$
(1+\varepsilon A)^{-1}X\to X\quad\mbox{\rm strongly in }\;L^p(\Omega\times [0,T]\times\mathcal O),
$$
see \cite{Aubin}.
Next by assumption \eqref{e1.2} we have
$$
\begin{array}{l}
\ds|\beta((1+\varepsilon A)^{-1}X)-\beta(X)|\\
\\
\ds \le\int_0^1\beta'(\lambda(1+\varepsilon A)^{-1}X)+(1-\lambda)
X)|(1+\varepsilon A)^{-1}X-X|d\lambda\\
\\
\le C\left(|(1+\varepsilon A)^{-1}X|^{m-1}+|X|^{m-1}+1\right) |(1+\varepsilon A)^{-1}X-X|.
\end{array}
$$
This yields, via H\"older's inequality  
$$
\begin{array}{l}
\ds\left|\E\int_0^t\int_\mathcal O(\beta((1+\varepsilon A)^{-1}X)-\beta(X))(X_\varepsilon-X)ds\, d\xi\right|\\
\\
\ds \le C|X_\varepsilon-X|_{L^p(\Omega\times [0,T]\times\mathcal O)}
|(1+\varepsilon A)^{-1}X-X|_{L^p(\Omega\times [0,T]\times\mathcal O)}
\\
\\
\ds \times \left(|(1+\varepsilon A)^{-1}X|^{m-1}_{L^p(\Omega\times [0,T]\times\mathcal O)}+
 |X|^{m-1}_{L^p(\Omega\times [0,T]\times\mathcal O)} +1 \right) \\
\\
\ds\le C_1|(1+\varepsilon A)^{-1}X-X|_{L^p(\Omega\times [0,T]\times\mathcal O)}\to 0,
\end{array}
$$
because $\{X_\varepsilon\}$ is bounded in $L^p(\Omega\times [0,T]\times\mathcal O)$ and $(m-1)\frac{p}{p-2}\le p$. Now the assertion follows by
\eqref{e3.24bis}.
\bigskip

Consider now the function
$$
\varphi(x)=\frac1p\;|x^-|^p_p.
$$
For any $x\in L^p(\mathcal O)$,  $\varphi$ is G\^ateaux differentiable and its differential $D\varphi\colon
L^p(\mathcal O)\to L^{p/(p-1)}(\mathcal O)$ is given by
$$
D\varphi(x)=-(x^-)^{p-1},
$$
while the second G\^ateaux derivative $D^2\varphi(x)\in L(L^p(\mathcal O);L^{p/(p-1)}(\mathcal O))$ is given by
$$
(D^2\varphi(x)h,g)=(p-1)\int_\mathcal O h\,g\;|x^-|^{p-2} d\xi,\quad \forall\; h,g,x\in L^p(\mathcal O).
$$
\begin{Lemma}
\label{l3.5}
Let $n\le 3$. For each $x\in L^p(\mathcal O)$ we have
\begin{equation}
\label{e3.24}
\begin{array}{l}
\ds\E[\varphi(X_\varepsilon(t))]+\E\int_0^t(A_\varepsilon X_\varepsilon(s),D\varphi(X_\varepsilon(s))ds
\\
\\
\ds=\varphi(x)+\frac{p-1}2\;\sum_{k=1}^\infty\mu^2_k\,\E\int_0^t\int_\mathcal O 
|X^-_\varepsilon(s)e_k|^2|X^-_\varepsilon(s)|^{p-2}dsd\xi.
\end{array}
\end{equation}

\end{Lemma}
{\bf Proof}. We note first that since $X_\varepsilon\in L_W^\infty(0,T;L^p(\Omega;L^p(\mathcal O)))$
the above formula makes sense. Next we approximate $\varphi$ by
$$
\varphi_\lambda(x)=\varphi((1+\lambda A_0)^{-1}x),\quad A_0=-\Delta,\;\; D(A_0)=
H^2(\mathcal O)\cap H_0^1(\mathcal O),\;\lambda>0.
$$
Since $\varphi\in C^2(C(\overline{\mathcal O}))$ and $(1+\lambda A_0)^{-1}$ is
linear continuous from $L^2(\mathcal O)$ to
$C(\overline{\mathcal O})$ (due to our assumption $n\le 3$) we  infer that $\phi_\lambda\in C^{2}(L^2(\mathcal O))$
and its first order and second order differentials
 are given,  respectively, by
$$
D\varphi_\lambda(x)=D\varphi((1+\lambda A_0)^{-1}x))(1+\lambda A_0)^{-1},
$$
$$
(D^2\varphi_\lambda(x)h,k)=(D^2\varphi((1+\lambda A_0)^{-1}x))((1+\lambda A_0)^{-1}h,(1+\lambda A_0)^{-1}k)$$
for  $h,k\in L^2(\mathcal O),x\in L^2(\mathcal O).$ Note that if $x\in L^p(\mathcal O)$, then
$$
D\varphi_\lambda(x)=-(1+\lambda A_0)^{-1}(((1+\lambda A_0)^{-1}x)^-)^{p-1}.
$$
So, for $\lambda\to 0$ we have $\varphi_\lambda(x)\to\varphi(x)$ and
$
 D\varphi_\lambda(x)\to D\varphi(x)\quad\mbox{\rm in}\;L^{p/(p-1)}(\mathcal O).
$
Next we write It\^o's formula for $\varphi_\lambda$ in the space $L^2(\mathcal O)$ which makes sense
by Lemma \ref{l3.2}.

We get
$$
\begin{array}{l}
\ds\E[\varphi_\lambda(X_\varepsilon(t))]+\E\int_0^t(A_\varepsilon(X_\varepsilon(s)),D\varphi_\lambda(X_\varepsilon(s)))ds
=\varphi_\lambda(x)
\\
\\
\ds+\frac{p-1}2\;\sum_{k=1}^\infty\mu^2_k\,\E\int_0^t\int_\mathcal O  |((1+\lambda
A_0)^{-1}(X_\varepsilon(s)e_k)|^2\;|((1+\lambda A_0)^{-1}X_\varepsilon(s))^-|^{p-2}d\xi\,ds.
\end{array}
$$
This yields
\begin{equation}
\label{e3.25}
\begin{array}{l}
\ds\E[\varphi_\lambda(X_\varepsilon(t))]-\E\int_0^t\int_\mathcal O(1+\lambda
A_0)^{-1}(A_\varepsilon(X_\varepsilon(s)))
(((1+\lambda A_0)^{-1}X_\varepsilon(s))^-)^{p-1}d\xi ds
\\
\\
\ds=\varphi_\lambda(x)\\
\\
\ds+\frac{p-1}2\;\sum_{k=1}^\infty\mu^2_k\,\E\int_0^t\int_\mathcal O  |((1+\lambda
A_0)^{-1}X_\varepsilon(s))^-|^{p-2} |(1+\lambda A_0)^{-1}(X_\varepsilon(s)e_k)|^2d\xi\,ds.
\end{array}
\end{equation}

 We know that for $\lambda\to 0$, $(1+\lambda A_0)^{-1}X_\varepsilon
(s)\to X_\varepsilon(s)$ strongly in $L^p(\mathcal O)$ a.e. in $\Omega\times(0,T) $ and
$$
|(1+\lambda A_0)^{-1}X_\varepsilon|_p
\le |X_\varepsilon|_p,\quad\mbox{\rm a.e. in}\;\Omega\times (0,T).
$$
Then by the Lebesgue dominated convergence theorem we have
\begin{equation}
\label{e3.26}
\lim_{\lambda\to 0}(1+\lambda A_0)^{-1}X_\varepsilon=
X_\varepsilon\quad\mbox{\rm strongly in}\;
L^p(\Omega\times (0,T)\times\mathcal O).
\end{equation}
Similarly, since $A_\varepsilon(X_\varepsilon)\in 
L^p(\Omega\times (0,T)\times\mathcal O)$ we have for $\lambda\to 0$
$$
(1+\lambda A_0)^{-1}(A_\varepsilon(X_\varepsilon))\to A_\varepsilon(X_\varepsilon)
,\quad\mbox{\rm strongly in}\;L^p(\Omega\times (0,T)\times\mathcal O).
$$
and  
$$
((1+\lambda A_0)^{-1} X_\varepsilon)^-\to X^-_\varepsilon
,\quad\mbox{\rm strongly in}\;L^{p}(\Omega\times (0,T)\times\mathcal O).
$$
This yields
\begin{equation}
\label{e3.27}
\begin{array}{l}
\ds\lim_{\lambda\to 0}\E\int_0^t\int_\mathcal O(1+\lambda A_0)^{-1}(A_\varepsilon(X_\varepsilon(s)))(((1+\lambda
A_0)^{-1}X_\varepsilon(s))^-)^{p-1}d\xi ds\\
\\
\ds=\int_0^t\int_\mathcal OA_\varepsilon(X_\varepsilon(s))(X^-_\varepsilon(s))^{p-1}d\xi ds.
\end{array}
\end{equation}
Then, if $x\in L^p(\mathcal O)$ letting $\lambda\to 0$ in \eqref{e3.25}
we get (since by Fatou's lemma $\E\varphi(X_\varepsilon(t))\le \liminf_{\lambda\to
0}\E\varphi_\lambda (X_\varepsilon(t)),\;\forall\;t\ge 0$)
$$
\begin{array}{l}
\ds\E[\varphi(X_\varepsilon(t))]-\E\int_0^t\int_\mathcal O
A_\varepsilon(X_\varepsilon(s))(X^-_\varepsilon(s))^{p-1} d\xi ds
\\
\\
\ds=\varphi (x)+\frac{p-1}2\;\sum_{k=1}^\infty\mu^2_k\,
\E\int_0^t\int_\mathcal O 
\ds|X_\varepsilon(s)e_k|^2\;|X^-_\varepsilon(s)|^{p-2}d\xi ds,
\end{array}
$$
 and so \eqref{e3.24} follows. $\Box$\bigskip

We have by   \eqref{e3.24} and the definition of $Y_\varepsilon$ that for  $x\in L^p(\mathcal O), x\ge 0$,
$$
\begin{array}{l}
\ds \E[\varphi(X_\varepsilon(t))]+\E\int_0^t\int_\mathcal O\Delta\beta (Y_\varepsilon(s))
(X^-_\varepsilon(s))^{p-1} ds d\xi
\\
\\
\ds=\frac{p-1}2\;\sum_{k=1}^\infty\mu^2_k\,
\E\int_0^t\int_\mathcal O |X^-_\varepsilon(s)e_k|^2\;|X^-_\varepsilon(s)|^{p-2}d\xi
ds\\
\\
\ds\le C\E\int_0^t|X^-_\varepsilon(s)|_p^p ds.
\end{array}
$$
(Recall that $A_\varepsilon(X_\varepsilon)=-\Delta\beta(Y_\varepsilon)$.)

We  therefore have, taking into account that  $\Delta\beta(Y_\varepsilon)=
\frac1\varepsilon(Y_\varepsilon-X_\varepsilon)$,  
\begin{equation}
\label{e3.28}
\begin{array}{l}
\ds \frac1p\;\E|X_\varepsilon^-(t)|^p_p+\frac1\varepsilon\;\E\int_0^t \int_\mathcal O
(Y_\varepsilon(s)-X_\varepsilon(s))(X^-_\varepsilon(s))^{p-1}d\xi ds\\
\\
\ds\le C\E\int_0^t|X^-_\varepsilon(s)|_p^p ds.
\end{array}
\end{equation}
We have
\begin{equation}
\label{e3.29}
|Y_\varepsilon^-(t)|^p_{p}\le \int_\mathcal OX_\varepsilon(t))(-Y^-_\varepsilon(t))^{p-1}d\xi,\quad \P\mbox{\rm-a.s.},
\end{equation}
analogously to deriving \eqref{e3.22}, for $x\in L^p(\mathcal O)$. To see this multiply \eqref{e3.21} by $g(y)$ where
$$
g(y)\colon=\frac{-(y^-)^{p-1}}{1+\lambda(y^-)^{p-2}},
$$
to get (after integration by parts) that
$$
\int_\mathcal O\frac{(y^-)^p}{1+\lambda(y^-)^{p-2}}\;d\xi+\varepsilon\int_\mathcal O\langle\nabla\beta(y), \nabla g(y)\rangle_{\R^n}\;d\xi
=\int_\mathcal O\frac{x^-(-y^-)^3}{1+\lambda(y^-)^2}\;d\xi.
$$
Note that $g$ as a composition of two decreasing Lipschitz functions is Lipschitz and decreasing.
So, we can apply Lemma \ref{l3.1} to obtain
$$
\int_\mathcal O\frac{(y^-)^4}{1+\lambda(y^-)^2}\;d\xi\le \int_\mathcal O\frac{x^-(-y^-)^3}{1+\lambda(y^-)^2}\;d\xi
$$
and \eqref{e3.29} follows by taking $\lambda\to \infty$.
By \eqref{e3.29} we have
$$
-|Y_\varepsilon^-(t)|^p_{p}\ge \int_\mathcal O(X_\varepsilon^+(t)-X^-_\varepsilon(t))(Y^-_\varepsilon(t))^{p-1}d\xi
\ge - \int_\mathcal OX_\varepsilon^-(t)(Y^-_\varepsilon(t))^{p-1}d\xi
$$
and therefore
$
|Y_\varepsilon^-(t)|^p_{p}\le   |X_\varepsilon^-(t)|_p\; |Y_\varepsilon^-(t)|_p^{p-1}.
$
Hence $|Y_\varepsilon^-(t)|_{p}\le |X_\varepsilon^-(t)|_{p}$ and so
$$
\int_\mathcal OY_\varepsilon^-(t)(X^-_\varepsilon(t))^{p-1}d\xi
\le
|X_\varepsilon^-(t)|^{p-1}_p\;|Y_\varepsilon^-(t)|_p\le |X_\varepsilon^-(t)|^p_p.
$$
Inserting the latter into\eqref{e3.28}  and taking into account that $Y_\varepsilon X^-_\varepsilon
\ge -Y^-_\varepsilon X^-_\varepsilon$
we see that $\E|X_\varepsilon^-(t)|^p_p=0$, a.e. $t\ge 0$ i.e,
$X_\varepsilon^-(t)=0$ a.e.  and therefore $X_\varepsilon(t)\ge 0$ a.e..  
Taking into account Lemma \ref{l3.4} we infer that $X\ge 0$. This completes the proof
in the case when $\beta$ is strictly monotone. $\Box$\bigskip

To treat the general  case of $\beta$ satisfying \eqref{e1.2} we shall associate to \eqref{e1.4} the equation
\begin{equation}
\label{e3.30}
\left\{\begin{array}{l}
dX^\lambda(t)+A^\lambda X^\lambda(t)=\sigma(X^\lambda(t))dW(t),\quad t\ge 0,\\
\\
X^\lambda(0)=x,
\end{array}\right.
\end{equation}
where
$$
A^\lambda(x)=-\Delta(\beta(x)+\lambda x),\quad \lambda>0
$$
and
$$
D(A^\lambda)=\{x\in H^{-1}(\mathcal O)\cap L^1(\mathcal O):\;\beta(x)+\lambda x\in H_0^{1}(\mathcal O)\}
$$
According to the first part of the proof, for each $x\in L^p(\mathcal O),x\ge 0$ and $\lambda>0$, equation
\eqref{e3.30} has a unique strong solution $X^\lambda$ which is  nonnegative a.e. on $\Omega\times (0,T)\times
\mathcal O$.

On the other hand, applying the It\^o formula from  \cite[Theorem I 3.2]{KR} to the equation
$$
d(X^\lambda(t)-X(t))+(A^\lambda X^\lambda(t)-AX(t))dt=(X^\lambda(t)-X(t))dW(t)
$$
where $X$ is the solution to \eqref{e1.1}, we get after  some calculations that
$$
\begin{array}{l}
\ds\frac12\;\E|X^\lambda(t)-X(t)|^2_{-1}+\lambda\E\int_0^t\langle X^\lambda(s),X^\lambda(s)-X(s)\rangle_{-1}  ds\\
\\
\ds\le\frac12\;\sum_{k=1}^\infty\mu_k^2\,\E\int_0^t|(X^\lambda(s)-X(s))e_k|^2_{-1}ds.
\end{array}
$$
This yields (see \eqref{e1.9}), since $$\langle X^\lambda(s),X^\lambda(s)-X(s)\rangle_{-1}\ge \langle
X(s),X^\lambda(s)-X(s)\rangle_{-1},$$
$$
\E|X^\lambda(t)-X(t)|^2_{-1}\le C\E\int_0^t|X^\lambda(s)-X(s)|^2_{-1}ds +\lambda^2\E\int_0^t |X(s)|_{-1}^2ds.
$$
Since $X\in C_W([0,T];L^2(\Omega,L^2(\mathcal O))$, we infer via Gronwall's lemma that
$$
\lim_{X^\lambda\to 0}X^\lambda=X \quad\mbox{\rm  in }\;C_W([0,T];L^2(\Omega,L^2(\mathcal O))
$$
and so $X\ge 0$ a.e. in $\Omega\times (0,T)\times \mathcal O$ as  claimed.

The final part of the assertion in Theorem \ref{t2.2} follows by the continuity of sample paths, since $L^p(\mathcal
O)$ is dense in $H^{-1}(\mathcal O)$ and the continuity
of solutions $X=X(t,x)$ with respect to the initial data $x$ (which
follows via It\^o's formula in the proof of uniqueness).
$\Box$

\section{Concluding remarks}

\hspace{5mm}  
Assumption $1\le n\le 3$ is unnecessarily strong and was taken for convenience only. As a matter of fact,
under suitable conditions of the form
\eqref{e1.8} we expect that Theorem \ref{t2.2} can be established  for any dimension $n$. This will be the subject 
of a forthcoming
paper.

  2) Theorem \ref{t2.2} and its proof remain valid for time--dependent nonlinear functions $\beta=\beta(t,x)$ where
$\beta$ is monotonically increasing in $x$, satisfies \eqref{e1.2} uniformly with respect to $t$ and
is continuous in $t$.

   3) One might speculate however that nonnegativity of $X(t,x)$ for $x\ge 0$ follows directly in
$H^{-1}(\mathcal O)$ by taking instead of $\varphi(x)=\frac1p|x^-|^p_p$ a suitable  $C^2$-function on $H^{-1}(\mathcal
O)$ which is zero on the cone of positive $x\in H^{-1}(\mathcal O)$ but so far we failed to find such a function.

\end{document}